\documentclass {article}
\usepackage{authblk}
\usepackage{etex}
\usepackage[utf8]{inputenc}
\usepackage[english]{babel}
\usepackage{amsmath}
\usepackage{amsthm}
\usepackage{amssymb}
\usepackage{sectsty}
\usepackage{titlesec}
\usepackage{color}
\usepackage[color,matrix,arrow]{xy}
\usepackage{amsgen}
\usepackage{amstext}
\usepackage{amsbsy}
\usepackage{amsopn}
\usepackage{amsfonts}
\usepackage{eepic}
\usepackage{graphicx}
\usepackage{epsf}
\usepackage{pstricks}
\usepackage{float}
\usepackage{enumerate}
\usepackage{eqnarray}
\usepackage{faktor}
\xyoption{all}

\usepackage{pgf}
\usepackage{tikz-cd}
\usetikzlibrary{automata, arrows.meta, positioning,decorations.pathmorphing}

\newcommand{\footrecall}[1]{%
} 
\usetikzlibrary{snakes,shapes,arrows,automata}

\newcommand{\Mon}{\mbox{\rm Mon}}

\titleformat*{\section}{\large\bfseries}
\titleformat*{\subsection}{\normalsize \bfseries}

\newcommand{\N}{\mathbb{N}}
\newcommand{\Z}{\mathbb{Z}}
\newcommand{\F}{\mathbb{F}}

\newcommand{\Q}{\mathbb{Q}}

\newcommand{\End}{\text{End}}

\newcommand{\Per}{\text{Per}}

\newcommand{\Orb}{\text{Orb}}

\newcommand{\Aut}{\text{Aut}}

\newcommand{\plog}{\phi\text{-}\log}
\newcommand{\pslog}{\psi\text{-}\log}
\newcommand{\tplog}{\tilde\phi\text{-}\log}
\newcommand{\tpslog}{\tilde\psi\text{-}\log}
\newcommand{\Img}{\text{Im}}

\newcommand{\sgn}{\text{sgn}}
\newcommand{\tsBrCP}{\text{2BrCP}}

\newcommand{\mc}{\mathcal}

\theoremstyle{definition}
\newtheorem{theorem}{Theorem}[section]
\newtheorem{corollary}[theorem]{Corollary}

\newtheorem{proposition}[theorem]{Proposition}

\newtheorem{lemma}[theorem]{Lemma}

\newtheorem{remark}[theorem]{Remark}

\begin{document}
 
 
\title{On the dynamics of endomorphisms of the direct product of two free groups}
\author{Andr\'e Carvalho\\

Centro de Matemática, Faculdade de Ciências da Universidade do Porto

  R. Campo Alegre s/n
  
  4169-007 Porto, Portugal

andrecruzcarvalho@gmail.com}
\maketitle

\begin{abstract}
We prove that Brinkmann's problems are decidable for endomorphisms of $F_n\times F_m$: given $(x,y),(z,w)\in F_n\times F_m$ and $\Phi\in \End(F_n\times F_m)$, it is decidable whether there is some $k\in \N$ such that $(x,y)\Phi^k=(z,w)$ (or $(x,y)\Phi^k\sim(z,w)$). We also prove decidability of a two-sided version of Brinkmann's conjugacy problem for injective endomorphisms which, from the work of Logan, yields a solution to the conjugacy problem in ascending HNN-extensions of $F_n\times F_m$. Finally, we study the dynamics of automorphisms of $F_n\times F_m$ at the infinity, proving that  that their dynamics at the infinity is asymptotically periodic, as occurs in the free and free-abelian times free cases.
\end{abstract}

\section{Introduction}
The dynamical study of automorphisms and endomorphisms of groups has been a topic of interest for many in the past 30 years. In this paper, we will consider two different questions about the dynamics of endomorphisms of $F_n\times F_m$, with $n,m\geq 2$,  the direct product of two nonabelian free groups. We remark that the study of endomorphisms of free-abelian times free groups and their dynamics is developed in \cite{[DV13],[Car22a], [CD24]}. 

The direct product of two free groups is known to be a good source of undecidability results, mostly due to Mihailova's construction \cite{[Mih58]}.  However, endomorphisms of these groups have been classified in \cite{[Car23b]}, and positive decidability results were obtained: the Whitehead problem was proven to be decidable, fixed and periodic subgroups were proven to be computable (in case they are finitely generated) and in \cite{[Car23d]} it is proved that the triviality of the intersection of the subgroup fixed by a monomorphism with the subgroup fixed by an endomorphism of $F_n\times F_m$ can be decided, showing that, in some sense, subgroups fixed by endomorphisms of $F_n\times F_m$ are quite special among arbitrary subgroups. Some dynamical results about the behaviour around \emph{infinite} fixed points are obtained in \cite{[Car23b]}.

In the first part, we will solve \emph{Brinkamnn's problems} for endomorphisms of $F_n\times F_m$. Brinkmann's (equality) problem on $G$, BrP($G$), consists on deciding, on input two elements $x,y\in G$ and an endomorphism $\phi\in \End(G)$, whether $y$ belongs to the $\phi$-orbit of $x$, i.e., whether there is some $k\in \N$ such that $x\phi^k=y$. Similarly, Brinkmann's conjugacy problem on $G$, BrCP($G$), consists on deciding whether  there is some $k\in \N$ such that $x\phi^k\sim y$. 
Brinkmann proved in \cite{[Bri10]} that these problems were decidable for automorphisms of the free group.  This turned out to be particularly important in proving decidability of the conjugacy problem for free-by-cyclic groups in \cite{[BMMV06]}. In \cite{[Log22]}, Logan solved several variations of this problem for general endomorphisms and used them to solve the conjugacy problem in ascending HNN-extensions of the free group, generalizing the work in \cite{[BMMV06]}. In fact, it is proved in \cite{[CD23]} that Logan's results imply decidability of Brinkmann's problems for endomorphisms (not necessarily injective) of the free group. Kannan and Lipton had already solved  in \cite{[KL86]} (Theorem \ref{kannanlipton})  the problem of deciding whether, given an $n\times n$ matrix $Q$ of rational numbers and two vectors of rational numbers $x,y\in \Q^n$, there is a natural number $i\in \N$ such that $xQ^i=y$, which is more general than  Brinkmann's Problem for free-abelian groups. In \cite{[Car23c]}, the link between Brinkmann's conjugacy problem and the conjugacy problem in cyclic extensions of the group was extended to \emph{generalized} versions of the problems and in \cite{[Car23]}, the author studied a \emph{quantification} of Brinkmann's problem in the context of virtually free groups.

We prove that both problems are decidable for endomorphisms of $F_n\times F_m$.

\newtheorem*{brp}{Theorem \ref{brp}}
\begin{brp}
There exists an algorithm taking as input integers $n,m>1$, $\Phi\in \End(F_n\times F_m)$ and two elements $(x,y),(z,w)\in F_n\times F_m$ that outputs $k\in \N$ such that $(x,y)\Phi^k=(z,w)$ if such a $k$ exists and outputs \texttt{NO} otherwise.
\end{brp}
\newtheorem*{brcp}{Theorem \ref{brcp}}
\begin{brcp}
There exists an algorithm taking as input integers $n,m>1$, $\Phi\in \End(F_n\times F_m)$ and two elements $(x,y),(z,w)\in F_n\times F_m$ that outputs $k\in \N$ such that $(x,y)\Phi^k\sim(z,w)$ if such a $k$ exists and outputs \texttt{NO} otherwise.
\end{brcp}

Logan proved the decidability of the conjugacy problem for ascending HNN-extensions of free groups by reducing this problem to a two-sided version of Brinkmann's problem and the twisted conjugacy problem for injective endomorphisms of $F_n$. We remark that
it is shown in \cite{[SS11]} that almost all one relator groups with $3$ or more generators are subgroups of ascending HNN-extensions of free groups and that the conjugacy problem is open for general one-relator groups.
We  solve this two-sided version of Brinkmann's conjugacy problem as well as the twisted conjugacy problem for injective endomorphisms of $F_n\times F_m$ and, using Logan's method, this yields decidability of the conjugacy problem for ascending HNN-extensions of $F_n\times F_m$. 
\newtheorem*{conjaschnn}{Corollary \ref{conjaschnn}}
\begin{conjaschnn}
The conjugacy problem is solvable for ascending HNN-extensions of $F_n\times F_m$, for $n,m>1$.
\end{conjaschnn}

Ascending HNN-extensions of free-abelian groups (which are direct products of free groups of rank 1) have been considered before (see, for example, \cite{[Val19]}) and include interesting classes of groups such as higher Baumslag-Solitar groups \cite{[Sto96],[SS17]}.

In the last part we study the \emph{dynamics at the infinity} of automorphisms of $F_n\times F_m$.

For automorphisms of free groups, infinite fixed points were discussed by Bestvina and Handel in  \cite{[BH92]} and Gaboriau, Jaeger, Levitt and Lustig in \cite{[GJLL98]}. The dynamics of free groups automorphisms is proved to be asymptotically periodic in \cite{[LL08]} and the same was obtained for free-abelian times free groups in \cite{[Car22a]}. In \cite{[CS09a]}, Cassaigne and Silva study the dynamics of infinite fixed points for monoids defined by special confluent rewriting systems (which contain free groups as a particular case). This was also achieved by Silva in \cite{[Sil13]} for virtually injective endomorphisms of virtually free groups and in  \cite{[Car22a]} and \cite{[Car23b]} this was done for free-abelian times free groups and the direct product of two free groups,  respectively.

We denote by $F_n$ the free group of rank $n$ and its alphabet by $X=\{x_1,\ldots,x_n\}$. Given two words $u$ and $v$ on a free group, we write $u\wedge v$ to denote the longest common prefix of $u$ and $v$. The prefix metric on a free group is defined by
$$d(u,v)=\begin{cases}
2^{-|u\wedge v|} \text{ if $u\neq v$}\\
0 \text{ otherwise}
\end{cases}.$$
The prefix metric on a free group is in fact an ultrametric and its completion $(\widehat{F}_n,\widehat{d})$ is a compact space which can be described as the set of all finite and infinite reduced words on the alphabet $X \cup X^{-1}$. We will denote by $\partial F_n$ the set consisting of only the infinite words and call it the \emph{boundary} of $F_n$.

We consider $F_n\times F_m$ endowed with the product metric given by taking the prefix metric in each factor.  This  is also an ultrametric and $\widehat {F_n \times F_m}$ is homeomorphic to $\widehat {F_n}\times \widehat {F_m}$ by uniqueness of the completion (Theorem 24.4 in \cite{[SW70]}).

When seen as a CAT(0) cube complex, or alternatively, as a median algebra, this coincides with the \emph{Roller compactification} (see \cite{[Rol98],[CFI16],[Fio20],[Fio21],[Bow22]}). Indeed, the Roller boundary and the Gromov boundary coincide in the free group and the behaviour of the Roller compactification when taking direct products is the same as the one of the completion of metric spaces, i.e., denoting by $\bar X$ the Roller compactification of $X$, we have that
$$\bar X=\bigcup_{i=1}^m  \bar X_1\times\ldots\times  \bar X_m.$$

It is well known by a general topology result \cite[Section XIV.6]{[Dug78]} that every uniformly continuous mapping $\varphi$ between metric spaces admits a unique continuous extension $\widehat\varphi$ to the completion. The converse is obviously true in this case: if a mapping between metric spaces admits a continuous extension to the completion, since the completion is compact, then the extension must be uniformly continuous, and so does the restriction to the original mapping. 
Uniformly continuous endomorphisms with respect to this metric are described in \cite{[Car23b]} and it turns out that every automorphism of $F_n\times F_m$ is indeed uniformly continuous. Characterizing and studying some properties of uniformly continuous endomorphisms has been done before for other classes of groups (see for example \cite{[CS09b],[Sil10],[Sil13],[AS16],[Car22a]}).

An interesting property of this metric (and so, of this boundary) is that the uniformly continuous endomorphisms of $F_n\times F_m$ for this metric $d$ are precisely the coarse-median preserving endomorphisms for the product coarse median obtained by taking the median operator induced by the coarse-median operators given by hyperbolicity of $F_n$ and $F_m$:

\newtheorem*{uccmp}{Theorem \ref{uccmp}}
\begin{uccmp}
Coarse-median preserving endomorphisms of $F_n\times F_m$ are precisely the uniformly continuous ones.
\end{uccmp}

Coarse-median preservation turns out to be a useful tool to obtain interesting properties of automorphisms (see \cite{[Fio21]}), including finiteness results on the fixed subgroup of an automorphism. We remark that uniformly continuous preserving endomorphisms also coincide with the uniformly continuous ones for free-abelian times free groups.

In the study of dynamical systems, the notion of \emph{$\omega$-limit set} plays a crucial role. Given a metric space $X$, a continuous function $f:X\to X$, and a point $x\in X$, the $\omega$-limit set $\omega(x,f)$ of $x$ consists of the accumulation points of the sequence of points in the orbit of $x$. Understanding the $\omega$-limits gives us a grasp on the behaviour of the system in the long term. If the space $X$ is compact, then $\omega$-limit sets are nonempty, compact and $f$-invariant.

In \cite{[LL08]}, the authors proved that in the case where $f$ is the extension of a free group automorphism to the completion, then, for every point $x\in \widehat F_n$, $\omega(x,f)$ is a periodic orbit.

In this paper, we prove the same result for automorphisms of $F_n\times F_m$:
\newtheorem*{recper}{Theorem \ref{recper}}
\begin{recper}
The continuous extension of every automorphism of $F_n\times F_m$ to the completion obtained by taking the prefix metric in each factor has asymptotically periodic dynamics.
\end{recper}

\section{Preliminaries}
The purpose of this section is to introduce the classification of endomorphisms of $F_n\times F_m$ obtained in  \cite{[Car23b]}, notation and some known results on the dynamics of endomorphisms of free and free-abelian groups that will be used later.

\subsection{Endomorphisms of $F_n\times F_m$} 

 For $(i,j)\in[n]\times [m]$, we define   $\lambda_i:F_n\to \mathbb Z$ as the homomorphism given by $a_k\mapsto \delta_{ik}$ and $\tau_j:F_m\to \mathbb Z$ given by $b_k\mapsto \delta_{jk}$, where $\delta_{ij}$ is the Kronecker symbol.
 
For $x\in F_n, y\in F_m$ and integers $p_i,q_i,r_j,s_j\in \Z$, we denote $\sum\limits_{i\in[n]} \lambda_i(x)p_i$ by $x^{P}$; $\sum\limits_{j\in[m]} \tau_j(y)r_j$ by $y^{R}$; $\sum\limits_{i\in[n]} \lambda_i(x)q_i$ by $x^{Q}$ and $\sum\limits_{j\in[n]} \tau_j(y)s_i$ by $y^{S}$. 
We also define $P=\{p_i\in\Z\mid i\in[n]\}$, $Q=\{q_i\in\Z\mid i\in[n]\}$, $R=\{r_j\in\Z\mid j\in[m]\}$ and $S=\{s_j\in\Z\mid j\in[m]\}$.
 We will keep this notation throughout the paper.

In \cite{[Car23b]}, the author classified endomorphisms of the direct product of two finitely generated free groups $F_n\times F_m$, with $m,n>1$ in seven different types:
\begin{enumerate}[(I)]
\item $(x,y)\mapsto\left(u^{x^P+y^R},v^{x^Q+y^S}\right)$, for some $1\neq u \in F_n$, $1\neq v \in F_m$ and integers $p_i,q_i, r_j, s_j\in \mathbb Z$ for $(i,j)\in [n]\times [m]$, such that $P,Q,R,S\neq \{0\}$. 
\item $(x,y)\mapsto\left(y\phi,v^{x^Q+y^S}\right)$, for some nontrivial homomorphism $\phi:F_m\to F_n$, $1\neq v \in F_m$  and integers $q_i,  s_j\in \mathbb Z$ for $(i,j)\in [n]\times [m]$,  such that $Q,S\neq \{0\}$. 
\item $(x,y)\mapsto\left(u^{x^P+y^R},y\phi\right),$  for some nontrivial endomorphism $\phi\in\End(F_m)$, $1\neq u \in F_n$, and integers $p_i,  r_j\in \mathbb Z$ for $(i,j)\in [n]\times [m]$,  such that $P,R\neq \{0\}$. 
\item $(x,y)\mapsto(y\phi,y\psi)$, for some nontrivial homomorphism $\phi: F_m\to F_n$ and nontrivial endomorphism $\psi\in \End(F_m)$.
\item $(x,y)\mapsto\left(1,v^{x^Q+y^S}\right)$, for some $1\neq v \in F_m$, and integers $q_i,  s_j\in \mathbb Z$ for $(i,j)\in [n]\times [m]$, such that $Q,S\neq \{0\}$.
\item $(x,y)\mapsto (x\phi,y\psi)$, for some endomorphisms $\phi\in \End(F_n)$, $\psi\in \End(F_m)$.
\item $(x,y)\mapsto (y\psi,x\phi)$, for homomorphisms $\phi:F_n\to F_m$  and $\psi:F_m\to F_n$.
\end{enumerate}

 From \cite[Proposition 3.2]{[Car23b]}, injective endomorphisms correspond to endomorphisms of type VI or VII such that the component mappings $\phi$ and $\psi$ are injective. We denote by $\Mon(G)$ the monoid of monomorphisms of a group $G$. In \cite[Corollary 3.3]{[Car23b]}, it is shown that automorphisms of $F_n\times F_m$ are the type VI endomorphisms with bijective component endomorphisms $\phi$ and $\psi$ and if $n=m$, there are also automorphisms of type VII, given by the type VII endomorphisms with bijective component homomorphisms $\phi$ and $\psi$. Hence, groups of the form $F_n\times F_n$ have \emph{more} automorphisms than groups of the form $F_n\times F_m$ with $n\neq m$. We will usually denote endomorphisms (and homomorphisms) of free groups by $\phi,\psi$ and endomorphisms of $F_n\times F_m$ by $\Phi,\Psi$.

Given an endomorphism by the image of the generators of $F_n\times F_m$, its type is decidable. Indeed, consider an endomorphism $\varphi: F_n\times F_m\to F_n\times F_m$ defined by $(a_i,1)\mapsto (x_i,y_i)$ and $(1,b_j)\mapsto (z_j,w_j)$ for $i\in [n]$ and $j\in [m].$ We define $X=\{x_i\mid i\in[n]\}$, $Y=\{y_i\mid i\in[n]\}$, $Z=\{z_j\mid j\in[m]\}$ and $W=\{w_j\mid j\in[m]\}$. We say that these sets are \emph{trivial} if they are singletons containing only the empty word and \emph{nontrivial} otherwise. As seen in \cite{[Car23b]}, matching the numbering above the endomorphisms are classified as follows:

\begin{enumerate}[(I)]
\item All sets $X,Y,Z$ and $W$ are nontrivial
\item $X$ is the only trivial set
\item $Y$ is the only trivial set
\item $X$ and $Y$ are the only trivial sets
\item $X$ and $Z$ are the only trivial sets
\item $Y$ and $Z$ are trivial sets
\item $X$ and $W$ are trivial sets
\end{enumerate}

Therefore, when we take an endomorphism as input (meaning that we are given images of the generators), we will often assume that its type is known.

\subsection{Dynamics of endomorphisms}
\subsubsection{Orbits}
We now review some concepts that will be used when studying the dynamics of endomorphisms.

Brinkmann's problems were proven to be decidable for automorphisms of the free group by Brinkmann in \cite{[Bri10]}. In \cite{[Log22]}, Logan extended these results to injective endomorphisms, and using a result from \cite{[Log22]}, it is shown in \cite{[CD23]} that these problems are decidable for arbitrary endomorphisms of the free group. Kannan and Lipton proved (something more general than) the decidability of Brinkmann's problem for free-abelian groups.
\begin{theorem}[Kannan--Lipton \cite{[KL86]}] \label{kannanlipton}
Given a matrix $M\in \mc M_{n\times n}(\Q)$ and two vectors $x,y\in \Q^n$ it is decidable whether there exists some $k\in \N$ such that $xM^k=y$.
 \end{theorem}
\begin{remark}\label{kannanlipton affine}
We will often use the above theorem for \emph{affine} transformations. It is seen in \cite{[CD24]} that this is not a problem, as an affine transformation can be seen as a restriction of a linear one. 
\end{remark}

Let $G$ be a group, $x,y \in G$ and $\phi \in \End (G)$. Then, the set of \emph{$\phi$-logarithms} (resp. \emph{$\widetilde{\phi}$-logarithms}) of $y$ in base $x$ is
$\plog_x(y)=\{k\geq 0\mid x\phi^k=y\}$ (resp. $\tplog_x(y)=\{k\geq 0\mid x\phi^k\sim y\}$). An element $x\in G$ for which there is some $k>0$ such that $x\phi^k=x$ (resp. $x\phi^k\sim x$) is said to be \emph{$\phi$-periodic} (resp. \emph{$\tilde\phi$-periodic}). The smallest $k$ satisfying the condition is called the \emph{$\phi$-period} (resp. \emph{$\tilde \phi$-period}).

Given a finite orbit $\Orb_\phi(x)$, we say that $\Orb_\phi(x) \cap\Per(\phi)$ is the \emph{periodic part of the orbit} and $\Orb_\phi(x)\setminus\Per(\phi)$ is the \emph{straight part of the orbit}. 

\begin{figure}[H]
\centering
  \begin{tikzcd}
    &&&&& \cdots \ar[ld,bend right] \\
    x\ar[r]
    & x\phi \ar[r]
    & x\phi^2 \ar[r]
    & \cdots  \ar[r]
    & x\phi^r \ar[rd, bend right] 
    && x\phi^{r+2} \ar[lu, bend right] \\
    &&&&& x\phi^{r+1} \ar[ru, bend right]
  \end{tikzcd}
  \caption{A finite orbit}
\end{figure}
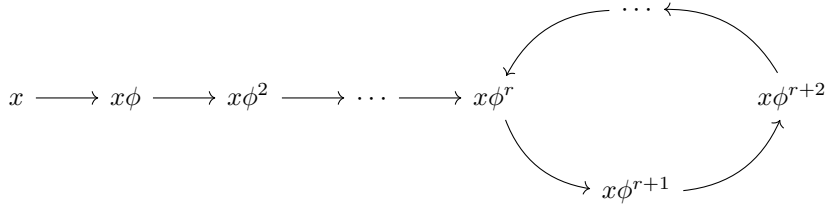
 
 In Figure 1, the \emph{straight part} of the orbit corresponds to $\{x,x\phi,\ldots,x\phi^{r-1}\}$ and the \emph{periodic part} of the orbit to $\{x\phi^k\mid k\geq r\}=\{x\phi^r,\ldots, x\phi^{r+p-1}\}$, where $p$ is the period of $x\phi^r$.
 We will consider the same notions up to conjugacy and the notation should be clear from context: for example, a \emph{finite $\tilde\phi$-orbit} is one containing only finitely many conjugacy classes. Naturally, $\phi$-logarithms and $\tilde\phi$-logarithms have the form $p+q\N$ for natural numbers $p$ and $q$ and are computable, as long as Brinkmann's problems are decidable.

\subsubsection{Dynamics at the infinity}

When we consider a hyperbolic group $G$ endowed with a visual metric $d$ and take its completion, we obtained the \emph{Gromov completion} of $G$. As said above, the endomorphisms of $G$ admitting a continuous extension to the completion are precisely the uniformly continuous ones (with respect to $d$). 

Now, let $f$ be a homeomorphism of a compact space $K$. Given $y\in K$, the \emph{$\omega$-limit set} $\omega(y,f)$, or simply $\omega(y)$, is the set of limit points of the sequence $f^n(y)$ as $n\to + \infty$. We say that the dynamics is \emph{asymptotically periodic} if every $\omega$-limit set is a periodic orbit (see \cite{[LL08]}).

We will study these concepts for $\widehat \Phi$, the continuous extension of an endomorphism to the completion $\widehat{F_n\times F_m}$ when the group is endowed with the product metric given by taking the prefix metric in each direct factor. As noted above, this is the \emph{Roller completion} of the group.

\begin{theorem}(Levitt--Lustig, \cite[Theorem I]{[LL08]})
\label{levittlustig}
Every automorphism of a free group has asymptotically periodic dynamics.
\end{theorem}

\section{Brinkmann's Problems}
The goal of this section is to prove the following theorems.
\begin{theorem}\label{brp}
There exists an algorithm taking as input integers $n,m>1$, $\Phi\in \End(F_n\times F_m)$ and two elements $(x,y),(z,w)\in F_n\times F_m$ that outputs $k\in \N$ such that $(x,y)\Phi^k=(z,w)$ if such a $k$ exists and outputs \texttt{NO} otherwise.
\end{theorem}
\begin{theorem}\label{brcp}
There exists an algorithm taking as input integers $n,m>1$, $\Phi\in \End(F_n\times F_m)$ and two elements $(x,y),(z,w)\in F_n\times F_m$ that outputs $k\in \N$ such that $(x,y)\Phi^k\sim(z,w)$ if such a $k$ exists and outputs \texttt{NO} otherwise.
\end{theorem}

We now present two technical lemmas that will be useful later.

\begin{lemma}\label{decide  conjugate power}
There exists an algorithm that, on input two reduced words $u,v\in F_n$, decides whether $u$ is conjugate to some power of $v$ and, in case it is, outputs the unique value of $k\in \N$ such that $u\sim v^k$.
\end{lemma}
\noindent\textit{Proof.} Two words $u$ and $v$ are conjugate if and only if the cyclic reduced core of $u$, $\tilde u$ is a cyclic permutation of the cyclic reduced core of $v$, $\tilde v$. In particular, their cyclically reduced cores must have the same length. So, we compute $\tilde u$,  $\tilde v$ and check if $k=\frac{|\tilde u|}{|\tilde v|}$ is an integer. If it is not, then there is no power $k$ such that $u\sim v^k$ and if it is, then it is our only candidate. Hence, it only remains to check whether $u\sim v^k$ or not.
\qed\\

The equality version of this lemma can be seen to hold in the same way.
\begin{lemma}\label{decide power}
There exists an algorithm that, on input two reduced words $u.,v\in F_n$, decides whether $u$ is equal to to some power of $v$ in $F_n$ and, in case it is, outputs the unique value of $k\in \N$ such that $u= v^k$.
\end{lemma}

Since the word problem and the conjugacy problem are decidable in $F_n\times F_m$ (as they reduce to the same problem in each factor), we will prove that we can decide the existence of a positive $k$ and that is equivalent to deciding if there is a nonnegative value of $k$ such that $(x,y)\Phi^k=(z,w)$ (or $(x,y)\Phi^k\sim(z,w)$, in the conjugacy case). 

\subsection{Type I} In this case, we have that $(x,y)\Phi=(u^{x^P+y^R},v^{x^Q+y^S})$, for some words $u,v\neq 1$.
    It is shown in \cite[Subsection 5.1]{[Car23b]} that, defining sequences in $\Z$ by 
\begin{align*}
\begin{cases}
a_1(x,y)=x^P+y^R\\
b_1(x,y)=x^Q+y^S\\
a_{n+1}(x,y)=a_n(x,y)u^P+b_n(x,y)v^R\\
b_{n+1}(x,y)=a_n(x,y) u^Q+b_n(x,y) v^S
\end{cases},
\end{align*}
we  have that, for every $k>0$, $(x,y)\Phi^k=(u^{a_k(x,y)},v^{b_k(x,y)})$. 
We want to decide, given $(x,y),(z,w)\in F_n\times F_m$, whether there is some $k\in \N$ such that $(x,y)\Phi^k=(z,w)$.  Using Lemma \ref{decide power} to decide if $z$ is a power of $u$ and $w$ is a power of $v$. If in any of the cases the answer is no, then we will not have a positive solution to our problem; if both answer yes, we compute exponents $(a,b)\in \N^2$ such that $(z,w)=(u^a, v^b)$.

So, for all $n\in \N$, $$\begin{bmatrix} a_{n} \\  b_n\end{bmatrix}=\begin{bmatrix}u^P & v^R\\u^Q  & v^S\end{bmatrix}^{n-1} \begin{bmatrix} x^P+y^R\\  x^Q+y^S\end{bmatrix},$$
hence, putting $M_\Phi=\begin{bmatrix}u^P & v^R\\u^Q  & v^S\end{bmatrix}$, the problem can now be translated as the problem of deciding whether there is some $k\in \N$ such that 
$$M_\Phi^k  \begin{bmatrix} x^P+y^R \\  x^Q+y^S\end{bmatrix}= \begin{bmatrix}a \\ b\end{bmatrix},$$
which can be done by Theorem \ref{kannanlipton}. 

In the conjugacy case, we start by checking if $k=0$ is a solution to our problem by solving the conjugacy problem. If not, using Lemma \ref{decide  conjugate power}, we check if $z$ is conjugate to a power of $u$ and $w$ is conjugate to a power of $v$.  If in any of the cases the answer is no, then we will not have a positive solution to our problem; if both answer yes, we compute the unique exponents $(a,b)\in \N^2$ such that $(z,w)\sim(u^a, v^b)$. Then we simply check if there is some $k$ such that $(x,y)\Phi^k=(u^a,v^b)$ using the equality algorithm above.

\subsection{Type II}
 In this case, we have that $(x,y)\Phi=(y\phi,v^{x^Q+y^S})$, for some word $v\neq 1$ and homomorphism $\phi:F_m\to F_n$.
It is shown in \cite[Subsection 5.2]{[Car23b]} that, defining a sequence in $\Z$ by 
\begin{align*}
\begin{cases}
a_1(x,y)=x^Q+ y^S\\
a_2(x,y)=a_1(x,y) v^S +y\phi^Q\\
a_n(x,y)=a_{n-1}(x,y)v^S +a_{n-2}(x,y)(v\phi)^Q\\
\end{cases}
\end{align*}
we  have that, for every $k>1$, $(x,y)\Phi^k=(v^{a_{k-1}(x,y)}\phi,v^{a_k(x,y)})$. We start by checking if $k=1$ is a solution. Then, by Lemma \ref{decide power} we test if $z$ is a power of $v\phi$ and if $w$ is a power of $v$. If one of them is not, then there is no solution $k$. If both are, we compute $a,b\in \N$ such that $(z,w)=(v^a\phi, v^b)$. Our problem now becomes finding $k>0$ such that $\begin{bmatrix} b  & a\end{bmatrix}=\begin{bmatrix} a_{k+1}(x,y) & a_{k}(x,y)\end{bmatrix}$. Similarly to the type I case, we use Kannan-Lipton's algorithm to decide the existence of such a $k$.
Consider the matrix
$M_\Phi=\begin{bmatrix}v^S & (v\phi)^Q\\1  & 0\end{bmatrix}$. Clearly, for $k>0$, 
$$M_\Phi\begin{bmatrix} a_{k+1}(x,y) \\ a_{k}(x,y)\end{bmatrix}=\begin{bmatrix} a_{k+2}(x,y) \\ a_{k+1}(x,y)\end{bmatrix},$$ and so 
$$M_\Phi^k\begin{bmatrix} a_{2}(x,y) \\ a_{1}(x,y)\end{bmatrix}=\begin{bmatrix} a_{k+2}(x,y) \\ a_{k+1}(x,y)\end{bmatrix}.$$
By Theorem \ref{kannanlipton}, we can decide if there is some $k$ such that   
$$\begin{bmatrix} a_{k+2}(x,y) \\ a_{k+1}(x,y)\end{bmatrix}=M_\Phi^k \begin{bmatrix} a_{2}(x,y) \\ a_{1}(x,y)\end{bmatrix}=\begin{bmatrix} b \\ a\end{bmatrix},$$ and we are done.

The variation up to conjugacy is similar to the previous case: we first decide if $z$ is conjugate to some power of $v\phi$ and $w$ is conjugate to some power of $v$ and then apply the algorithm for equality.
\subsection{Type III}
We have that $(x,y)\mapsto\left(u^{x^P+y^R},y\phi\right)$,  for some word $u\neq 1$ and endomorphism $\phi\in\End(F_m)$. It is seen in \cite[Subsection 5.3]{[Car23b]} that  $$(u^a,y)\Phi^k=\left(u^{a(u^P)^{k}+\sum\limits_{t=0}^{k-1} (y\phi^t)^R(u^P)^{k-t-1}},y\phi^k\right).$$
We start by solving BrP($F_n$) to compute $p_0, p_1\in \N$ such that $\plog_y(w)=p_0+p_1\N$.  
If $p_1=0$, then $p_0$ is our only candidate, so the only thing to check is if $(x,y)\Phi^{p_0}=(z,w)$. If $p_1\neq 0$, $w$ is a $\phi$-periodic point with period $p_1$. We check if $z$ is a power of $u$. If it is not, then we are done, since there is no possible $k$; if it is, we compute $a\in \N$ such that $z=u^a$. 

Now consider the following affine transformation of $\Z$
\begin{align*}
\Theta\colon \Z&\longrightarrow\Z\\
c&\longmapsto  c(u^P)^{p_1}+\sum_{t=0}^{p_1-1} (w\phi^{t})^R(u^P)^{p_1-t-1},
\end{align*}
and denote by $(a_n)_n$ the sequence such that $(x,y)\Phi^n=(u^{a_n},y\phi^n)$. Clearly, given $k\in \N$, we can compute $a_k$. We claim that, for $k\in p_0+p_1\N$, writing $k=p_0+p_1r$, we have that  
\begin{align}\label{claimIII}
(x,y)\Phi^k=(u^{a_{p_0}\Theta^r}, y\phi^k).
\end{align}
We prove it by induction over $r$. Clearly, if $r=0$, then $(x,y)\Phi^k=(x,y)\Phi^{p_0}=(u^{a_{p_0}}, y\phi^k).$ Now assume that the claim holds for $r$ up to some $n\geq 0$. Notice that $w=y\phi^{p_0}$ is periodic with period $p_1$, and so, for all $n, t\in \N$, 
$y\phi^{p_0+np_1+t}=w\phi^t.$ Hence, 
\begin{align*}
(x,y)\Phi^{p_0+(n+1)p_1}&=(x,y)\Phi^{p_0+np_1}\Phi^{p_1}\\
&=(u^{a_{p_0}\Theta^n}, y\phi^{p_0+np_1})\Phi^{p_1}\\
&=\left(u^{(a_{p_0}\Theta^n)(u^P)^{p_1}+\sum\limits_{t=0}^{p_1-1} (y\phi^{p_0+np_1+t})^R(u^P)^{p_1-t-1}},y\phi^{p_0+(n+1)p_1}\right)\\
&=\left(u^{(a_{p_0}\Theta^n)(u^P)^{p_1}+\sum\limits_{t=0}^{p_1-1} (w\phi^{t})^R(u^P)^{p_1-t-1}},y\phi^{p_0+(n+1)p_1}\right)\\
&=\left(u^{a_{p_0}\Theta^{n+1}}, y\phi^k\right).
\end{align*}
Since our only candidate solutions are of the form $p_0+p_1\N$, we can proceed as follows: we check manually the existence of a solution for $k$ up to $p_0$, using the word problem; then, if we obtain no positive answer, we verify the existence of some $r\in \N$ such that $a_{p_0}\Theta^r=a$, which can be done by Remark \ref{kannanlipton affine}.

When considering the conjugation variant, we proceed analogously. Using BrCP($F_n$), we compute  $p_0, p_1\in \N$ such that $\tplog_y(w)=p_0+p_1\N$. We then check if $z$ is conjugate to some power of $u$ and if so, compute $a$ such that $u^a\sim z$. Using the same argument as above, we check if there is some $r$ such that the first component $(x,y)\Phi^{p_0+rp_1}$ is $u^a$. In fact, even though $w$ is no longer a periodic point, it is enough that $w$ is $\tilde\phi$-periodic since, for words $u,v\in F_m$ such that $u\sim v$, i.e., for which there is some $z\in F_m$ such that $u=z^{-1}vz$, we have that $$u^R=\sum\limits_{j\in[m]} \tau_j(u)r_j=\sum\limits_{j\in[m]} \tau_j(z^{-1}vz)r_j=\sum\limits_{j\in[m]} \tau_j(v)r_j=v^R.$$ This means that for all $n, t\in \N$, 
$(y\phi^{p_0+np_1+t})^R=(w\phi^t)^R$ and the same reduction done above can be done here and the claim (\ref{claimIII}) holds in this case too. Hence, we reduce our problem to finding some $r\in \N$ such that $a_{p_0}\Theta^r=a$, which can be done by Remark \ref{kannanlipton affine}.

\subsection{Type IV}
We have that  $(x,y)\mapsto(y\phi,y\psi)$, for some homomorphism $\phi: F_m\to F_n$ and endomorphism $\psi\in \End(F_m)$, and so, as seen in \cite[Subsection 5.4]{[Car23b]}, $(x,y)\Phi^k=(y\psi^{k-1}\phi, y\psi^k),$ for $k>0$. Using BrP($F_n$), we start by computing $\pslog_y(w)$. If it is empty, then there are no solutions $k$. If not,  $\pslog_y(w)=p_0+p_1\N$ for some computable $p_0,p_1\in \N$. We check if $k=p_0$ is a solution or not. We have that $w=y\psi^{p_0}$ is a $\psi$-periodic point with period $p_1$. For $k=p_0+p_1r$, with $r\geq 1$, we have that  $(x,y)\Phi^k=(y\psi^{p_0+p_1r-1}\phi, w)$. But for all $r\geq 1$, $y\psi^{p_0+p_1r-1}=y\psi^{p_0+p_1-1}$. So we simply check if $z=y\psi^{p_0+p_1-1}\phi$. If it is, then $p_0+p_1\N_{>0}$ is the set of solutions; if not, then there are no solutions to our problem.

We can handle the conjugacy version of the problem in the same way, computing $\tpslog_y(w)$ instead of $\pslog_y(w)$. In this case, we have that,  for all $r\geq 1$, $y\psi^{p_0+p_1r-1}\sim y\psi^{p_0+p_1-1}$ and so $y\psi^{p_0+p_1r-1}\phi\sim y\psi^{p_0+p_1-1}\phi$  and the set of solutions is $p_0+p_1\N_{>0}$ if $z\sim y\psi^{p_0+p_1-1}\phi$ and empty otherwise.
\subsection{Type V}
We have that  $(x,y)\mapsto\left(1,v^{x^Q+y^S}\right)$, for some $1\neq v \in F_m$, for some word $v\in F_m$, and so, as seen in \cite[Subsection 5.5]{[Car23b]}, $(x,y)\varphi^k=\left(1,v^{(x^Q+y^S)(v^S)^{k-1}}\right)$, for $k>0$. 
So we start by checking if $z=1$ and if $w$ is a power of $v$. If any of these does not hold, then there is no positive solution $k$. If both hold, we compute $a\in \Z$ such that $w=v^a$. Clearly, $k$ is a solution if and only if $(v^S)^{k-1}=\frac{a}{x^Q+y^S}$, and its existence is decidable. The conjugacy problem is analogous: we compute the unique $a'\in \Z$ such that $w\sim v^{a'}$ and check if there is some $k$ such that $(v^S)^{k-1}=\frac{a'}{x^Q+y^S}$.
\subsection{Type VI}
We have that  $(x,y)\mapsto (x\phi,y\psi)$, for some endomorphisms $\phi\in \End(F_n)$, $\psi\in \End(F_m)$, and so t $(x,y)\varphi^k=(x\phi^{k},y\psi^k)$, for $k>0$. So this case is also simple: we compute $p_0,p_1, q_0,q_1\in \N$ such that $\plog_x(z)=p_0+p_1\N$ and $\psi\text{-}\log_y(w)=q_0+q_1\N$ (if there are no $p_0$ or $q_0$, then there are no solutions to our problem). Our set of solutions is the intersection $\plog_x(z)\cap \psi\text{-}\log_y(w)$, which we can check if it is empty or not. For the conjugacy version of the problem we proceed in the same way, computing $\tplog_x(z)$ and $\tilde\psi\text{-}\log_y(w)$ instead.

\subsection{Type VII}
We have that $(x,y)\mapsto (y\psi,x\phi)$, for homomorphisms $\phi:F_n\to F_m$  and $\psi:F_m\to F_n$, and so, as seen in  \cite[Subsection 5.7]{[Car23b]}, for $k\geq 0$,  $(x,y)\varphi^{2k}= (x(\phi\psi)^k,y(\psi\phi)^k)$ and $(x,y)\varphi^{2k+1}= (y(\psi\phi)^k\psi,x(\phi\psi)^k\phi)$. So we will describe two algorithms: one to decide the existence of an even $k$, and another one to decide the existence of an odd $k$. The first one is simple, since it can be seen as an application of the previous case, defining $\Phi$ to be the type VI endomorphism mapping $(x,y)$ to $(x\phi\psi, y\psi\phi)$. Now we describe the second one. We start by checking if there is some even $r\in \N$ such that $(x,y)\Phi^r=(z,w)\Phi$. If there is none, then there is no solution $k$; otherwise compute such a minimal $r$. It is proved in  \cite[Subsection 5.7]{[Car23b]} that the subgroup of periodic points of $\Phi$ is given by $\Per(\Phi)=\Per(\phi\psi)\times \Per(\psi\phi)$. In particular, it is computable and  we can check if a given element of $F_n\times F_m$ is periodic or not. We verify if $(w\psi,z\phi)=(z,w)\Phi$ is periodic or not. If it is not, then our only candidate is $k=r-1$. Indeed, if there was some other $k>r$, then we would have that $$\left((z,w)\Phi\right)\Phi^{k-r+1}=(x,y)\Phi^r\Phi^{k-r+1}=(x,y)\Phi^{k+1}=(z,w)\Phi$$ and if there was some other $k<r$ we would have that 
$$(z,w)\Phi\Phi^{r-k}=(x,y)\Phi^{k}\Phi^{r-k+1}=(x,y)\Phi^{r+1}=(z,w)\Phi.$$
So, assume that $(z,w)\Phi$ is periodic with period $s$. In this case, there are only finitely many checks to do since the orbit of $(x,y)$ is finite, as illustrated in the following picture:

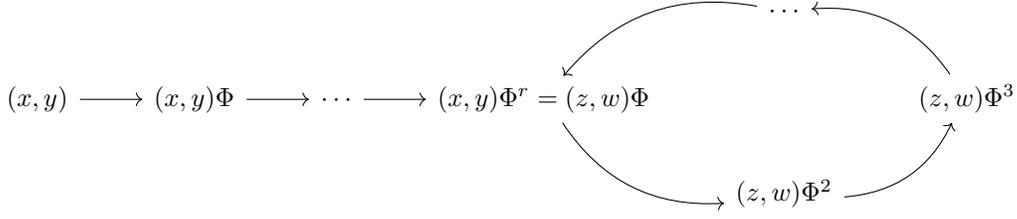
\begin{figure}[H]
\centering
  \begin{tikzcd}
    &&&& \cdots \ar[ld,bend right] \\
    (x,y)\ar[r]
    & (x,y)\Phi \ar[r]
    & \cdots  \ar[r]
    &(x,y)\Phi^r=(z,w)\Phi\ar[rd, bend right] 
    && (z,w)\Phi^3 \ar[lu, bend right] \\
    &&&& (z,w)\Phi^2 \ar[ru, bend right]
  \end{tikzcd}
  \caption{The orbit of $(x,y)$.}
\end{figure}

The conjugacy case is similar. Checking the existence of an even $k$ can be translated as an instance of the algorithm for Type VI endomorphisms in the same way as above. Then we verify if there is some even $r$ such that $(x,y)\Phi^r\sim (z,w)\Phi$. If not, then there is no solution $k$, and if there is one, we compute it. We can check if $(z,w)\Phi$ is $\tilde\Phi$-periodic, because that happens if and only if there is some $s\in 2\N$ such that $(z,w)\Phi\Phi^2\Phi^s\sim (z,w)\Phi$. Indeed, if there is such an $s$, then $(z,w)\Phi\Phi^{s+2}\sim (z,w)\Phi$ and $(z,w)\Phi$ is $\tilde \Phi$-periodic and conversely, if $(z,w)\Phi$ is $\tilde\Phi$-periodic, then there is some $p\in \N_{>0}$ such that $(z,w)\Phi\Phi^p\sim (z,w)\Phi$ and so
 $$(z,w)\Phi\Phi^2\Phi^{2p-2}=(z,w)\Phi\Phi^{2p}\sim (z,w)\Phi\Phi^p\sim (z,w)\Phi.$$ As done in the equality case, if $(z,w)\Phi$ is not $\tilde\Phi$-periodic, there is only one possible solution and if it is, our problem amounts to finitely many checks as the orbit of $(x,y)$ has only finitely many conjugacy classes.

\section{The conjugacy problem for ascending HNN-extensions of $F_n\times F_m$}

Let $G$ be a group and $\Phi$ be an injective endomorphism of $G$. The \emph{ascending HNN-extension of $G$ induced by $\Phi$} is the group with relative presentation 

\begin{equation} \label{eq: presentation AHNN}
    G \ast_\Phi
    \,=\,
    \langle
    G,
    x
\mid
    x^{-1}gx=\Phi(g), g\in G
    \rangle
\end{equation}
and with normal forms given by $x^i g x^{-j}$, with $i,j\in \N_0$ and $g\in G.$

In  \cite{[BMMV06]}, it was proved that the conjugacy problem was solvable in free-by-cyclic groups, by reducing this problem to the solution of BrCP($F_n$) and TCP($F_n$) for automorphisms.  In \cite{[Log22]}, Logan proved that the conjugacy problem was decidable for ascending HNN-extensions of a free group, extending the work in \cite{[BMMV06]}. To do so, Logan solved a variation of BrCP, which we call \emph{two-sided Brinkmann's conjugacy problem} and TCP($F_n$) for injective endomorphisms of the free group.
The \emph{two-sided Brinkmann conjugacy problem} in $G$ with respect to $\mathcal{T}\subseteq \End(G)$, denoted by $\tsBrCP(G)$ is the problem of deciding, on input an endomorphism $\Phi\in \mc T$ and elements $x,y\in G$, whether there is some pair $(r,s)\in \N^2$ such that $x\Phi^r \sim y\Phi^s$. 

We will show that TCP($F_n\times F_m$) and $\tsBrCP(F_n\times F_m)$ are decidable with respect to monomorphisms, which, using Logan's techniques implies that the conjugacy problem is decidable for ascending HNN-extensions of the direct product of two free groups. Notice that the free-abelian times free case was already dealt with in \cite{[CD24]}. We recall that it is proved in \cite[Proposition 3.2]{[Car23b]} that injective endomorphisms of $F_n\times F_m$ are precisely the types VI and VII endomorphisms having injective component homomorphisms $\phi$ and $\psi$.
We start by proving a lemma that will be useful in proving decidability of $\tsBrCP(F_n\times F_m)$ in the injective case.

\begin{lemma}\label{eventually phi periodic}
There exists an algorithm taking as input an integer $n>1$, an injective endomorphism $\phi\in \Mon(F_n)$, and an element $x\in F_n$ that decides if the number of conjugacy classes in the orbit of $x$ is finite or not. Equivalently, it is decidable whether there is a pair $(p,q)\in \N^2$ such that $x\phi^p\sim x\phi^q$ with $p>q$.
\end{lemma}
\noindent\textit{Proof.} First, we verify if $x$ is $\tilde\phi$-periodic using  \cite[Lemma 4.1]{[Log22]} to decide if there is some $k\in \N$ such that $x\phi\phi^k\sim x$. If $x$ is $\tilde\phi$-periodic, then the orbit of $x$ contains only finitely many conjugacy classes and we are done. So, we may assume that $x$ is not $\tilde\phi$-periodic.

Suppose that there are  $(p,q)\in \N^2$ such that $x\phi^p\sim x\phi^q$ with $p>q$ and let $p$ be the minimal positive integer for which there is some $q < p$ such that $x\phi^p\sim x\phi^q$. Then, there is some $u\in \F_n$ such that 
$x\phi^p=u^{-1}x\phi^qu$. Notice that $q>0$, since otherwise we would have that $x$ is $\tilde\phi$-periodic. If $u\in \Img(\phi)$, i.e., $u=u'\phi$ for some $u'\in F_n$, then $x\phi^p=(u'^{-1}\phi)(x\phi^q)(u'\phi)$, which, by injectivity implies that $x\phi^{p-1}=u'^{-1}x\phi^{q-1}u'$ and this contradicts the minimality of $p$. Thus, for the minimal $p$, we know that $x\phi^p=u^{-1} x\phi^q u$, for some $u\in F_n\setminus \Img(\phi)$. 

Therefore, the number of conjugacy classes in the orbit of $x$ is finite if and only if there is some pair $(p,q)\in \N^2$ such that $p>q>0$ and some $u\in F_n\setminus \Img(\phi)$ such that  $x\phi^p=u^{-1}x\phi^qu$. Verifying this condition can be done by applying \cite[Lemma 4.3]{[Log22]} to verify if there is some pair $(p,q)\in \N^2$ with $p\geq q \geq 0$ such that there is some element $u\in F_n\setminus \Img(\phi)$ satisfying $(x\phi)\phi^p=u^{-1}(x\phi^q)u$.
\qed\\
\begin{theorem}
There exists an algorithm taking as input two integers $n,m>1$,  two elements $(x,y),(z,w)\in F_n\times F_m$ and an injective endomorphism $\Phi\in \Mon(F_n\times F_m)$ that outputs a pair $(r,s)\in \N^2$ such that $(x,y)\Phi^r\sim(z,w)\Phi^s$ if such a pair exists and outputs \texttt{NO} otherwise.
\end{theorem}
\noindent\textit{Proof.} Suppose that $\Phi$ is an injective endomorphism of type VI. Then, $(x,y)\Phi=(x\phi,y\psi)$ for some $\phi\in \Mon(F_n)$ and $\psi\in \Mon(F_m)$. Naturally, 
$$(x,y)\Phi^r\sim(z,w)\Phi^s\iff (x\phi^r,y\psi^r)\sim(z\phi^s,w\psi^s)\iff \begin{cases} x\phi^r\sim z\phi^s\\ y\psi^r\sim w\psi^s \end{cases}.$$
We start by using \cite[Proposition 4.4]{[Log22]} to check whether there is a pair $(r_1,s_1)$ such that $x\phi^{r_1}\sim z\phi^{s_1}$. If there is none, then we answer \texttt{NO}. If there is such a pair, we compute the minimal $r_1$ for which there is some $k\in \N$ such that $x\phi^{r_1}\sim z\phi^{k}$. This is possible by \cite[Lemma 4.1]{[Log22]}: we keep increasing the candidate values $r_1$ until the algorithm in \cite[Lemma 4.1]{[Log22]} answers \texttt{YES}. When it does, we compute the minimal $s_1$ such that $x\phi^{r_1}\sim z\phi^{s_1}$, which can be done by iteratively solving the conjugacy problem.
 Then we verify if the orbit of $z$ has finitely many conjugacy classes or not, using Lemma \ref{eventually phi periodic}. If there are only finitely many conjugacy classes, we compute the minimal $q_z\in \N$ such that $z\phi^{q_z}$ is $\tilde\phi$-periodic and $p_z\in \N$, the $\tilde\phi$-period of $z\phi^{q_z}$. If there are infinitely many conjugacy classes, we say that $q_z=\infty$.
 We claim that the set $\mc S$ of all pairs $(p,q)\in \N^2$ such that $x\phi^{p}\sim z\phi^{q}$ is precisely
 \small{
$$\begin{cases}
(r_1,s_1)+(1,1)\N  &\text{ if $q_z=\infty$}\\
(r_1,s_1)+(1,1)\N+(0,p_z)\N + (p_z,0)\N &\text{ if $q_z\leq s_1$}\\
\{(r_1+i,s_1+i)\mid 0\leq i< q_z-s_1\}\cup \left((r_1+q_z-s_1,q_z)+(1,1)\N+(0,p_z)\N + (p_z,0)\N\right)&\text{ if $q_z> s_1$}
 \end{cases}$$}
 \normalsize
Suppose that $q_z=\infty$. Then, clearly $(r_1,s_1)+(1,1)\N\subseteq \mc S$. So, suppose that  $x\phi^{p}\sim z\phi^{q}$. Then, by minimality of $r_1$, we have that $p\geq r_1$. Clearly, 
\begin{align}\label{zqzs1pr1}
z\phi^q\sim x\phi^p=x\phi^{r_1}\phi^{p-r_1}\sim z \phi^{s_1}\phi^{p-r_1}=z\phi^{s_1+p-r_1}.
\end{align}
Since the orbit of $z$ has infinitely many conjugacy classes, then it must be the case that $q=s_1+p-r_1$, i.e., $q-s_1=p-r_1$, and so, $(p,q)=(r_1,s_1)+ (p-r_1)(1,1)\in (r_1,s_1)+ (1,1)\N$.

Now assume that $q_z\leq s_1$. Then $z\phi^{s_1}$ has $\tilde\phi$-period $p_z$, and so $x\phi^{r_1}$ must also have  $\tilde\phi$-period $p_z$. Indeed, $x\phi^{r_1}\phi^k\sim x\phi^{r_1}$ if and only if $z\phi^{s_1}\phi^k\sim  x\phi^{r_1}\sim z\phi^{s_1}$.
Thus, it is clear that $(r_1,s_1)+(1,1)\N+(0,p_z)\N+ (p_z,0)\N\subseteq \mc S$.
Now suppose that  $x\phi^{p}\sim z\phi^{q}$. Again, by minimality of $r_1$, we have that $p\geq r_1$, which implies that  $x\phi^{p}$ (and so, $z\phi^q$) is in the $\tilde\phi$-periodic part of the orbit, and has period $p_z$.
Proceeding as above, we get that $z\phi^q\sim z\phi^{s_1+p-r_1}$. If $q>s_1+p-r_1$, there is some $r\in \N$ such that $q-p-s_1+r_1=rp_z$, which means that $q=s_1+(p-r_1)+rp_z$, thus $$(p,q)=(r_1,s_1)+(p-r_1)(1,1)+ r(0,p_z)\in (r_1,s_1)+(1,1)\N+(0,p_z)\N+ (p_z,0)\N.$$

If, on the other hand $q\leq s_1+p-r_1$, there is some $r\in \N$ such that $s_1+p-r_1-q=rp_z$, and so $p=r_1+(q-s_1)+rp_z$, thus $$(p,q)=(r_1,s_1)+(q-s_1)(1,1)+ r(p_z,0)\in (r_1,s_1)+(1,1)\N+(0,p_z)\N+ (p_z,0)\N.$$

Now, suppose that $q_z>s_1$. Obviously $\{(r_1+i,s_1+i)\mid 0\leq i< q_z-s_1\}\subseteq \mc S$. Again, since $x\phi^{r_1+q_z-s_1}\sim z\phi^{q_z}$, then $x\phi^{r_1+q_z-s_1}$ has $\tilde\phi$-period $p_z$, and so it is clear that $(r_1+q_z-s_1,q_z)+(1,1)\N+(0,p_z)\N + (p_z,0)\N\subseteq \mc S$. Now suppose that  $x\phi^{p}\sim z\phi^{q}$. By minimality of $r_1$, we must have that $p\geq r_1$. Suppose that $p< r_1+q_z-s_1$. Then $s_1+p-r_1<q_z$, which means that $ z\phi^{s_1+p-r_1}$ is not $\tilde\phi$-periodic. 
 As in (\ref{zqzs1pr1}), $z\phi^q\sim z\phi^{s_1+p-r_1}$ and so, it must be that $q= s_1+p-r_1$ and $(p,q)=(r_1,s_1)+(p-r_1,q-s_1)=(r_1,s_1)+(p-r_1,p-r_1)$. Since we are assuming that $p-r_1<q_z-s_1$, it follows that $(p,q)\in \{(r_1+i,s_1+i)\mid 0\leq i< q_z-s_1\}.$
If, on the other hand,  $p \geq r_1+q_z-s_1$, $$x\phi^p\sim x\phi^{r_1}\phi^{p-r_1}\sim z\phi^{s_1}\phi^{p-r_1},$$
with $p-r_1\geq q_z-s_1$ and so $x\phi^p$  (and so, $z\phi^q$) lies in the $\tilde\phi$-periodic part of the orbit, and has period $p_z$. So we proceed as in the previous case. As in (\ref{zqzs1pr1}), $z\phi^q\sim z\phi^{s_1+p-r_1}$.
If $q>s_1+p-r_1$, there is some $r\in \N$ such that $q-p-s_1+r_1=rp_z$, which means that $q=s_1+(p-r_1)+rp_z$, thus $$(p,q)=(r_1+q_z-s_1,q_z)+(p-r_1+s_1-q_z)(1,1)+ r(0,p_z)\in (r_1,s_1)+(1,1)\N+(0,p_z)\N+ (p_z,0)\N.$$

If, on the other hand $q\leq s_1+p-r_1$, there is some $r\in \N$ such that $s_1+p-r_1-q=rp_z$, and so $p=r_1+(q-s_1)+rp_z$, thus $$(p,q)=(r_1+q_z-s_1,q_z)+(q-q_z)(1,1)+ r(p_z,0)\in (r_1,s_1)+(1,1)\N+(0,p_z)\N+ (p_z,0)\N.$$

Proceeding analogously, we compute the set $\mc S'\subseteq \N^2$ of pairs $(p,q)\in \N^2$ such that $y\phi^{p}\sim w\phi^{q}$. These are semilinear sets, so their intersection can be computed (in particular its emptyness can be decided) and the solution set for our problem is precisely $\mc S\cap \mc S'$.

Now suppose that $\Phi$ is a type VII injective endomorphism. Then $(x,y)\Phi=(y\psi,x\phi)$ for some injective homomorphisms $\phi:F_n\to F_m$ and $\psi:F_m\to F_n$ and, for $k\geq 0$,  $(x,y)\varphi^{2k}= (x(\phi\psi)^k,y(\psi\phi)^k)$ and $(x,y)\varphi^{2k+1}= (y(\psi\phi)^k\psi,x(\phi\psi)^k\phi)$. We prove that we can decide if there is a pair $(r,s)\in \N^2$ by showing that it is decidable whether there is such a pair $(r,s)$ in $2\N\times 2\N$, in $2\N\times (2\N+1)$, in $(2\N+1)\times 2\N$ and in $(2\N+1)\times (2\N+1)$.

The first case is immediate, since we are deciding if there are some $r',s'\in \N^2$ such that $(x(\phi\psi)^{r'},y(\psi\phi)^{r'})\sim (z(\phi\psi)^{s'},w(\psi\phi)^{s'})$, which can be reduced to the type VI case, by defining the type VI monomorphism given by $\Psi: (x,y)\mapsto (x\phi\psi,y\psi\phi)$. The second and third cases are the same since one can be obtained from the other by swapping the roles of $(x,y)$ and $(z,w)$. We will now see that we can decide whether there are some $r',s'\in \N$ such that $(x,y)\Phi^{2r'+1}\sim (z,w)\Phi^{2s'}$.
But 
\begin{align*}
&(x,y)\Phi^{2r'+1}\sim (z,w)\Phi^{2s'}\\
\iff\; &  (y(\psi\phi)^{r'}\psi,x(\phi\psi)^{r'}\phi)\sim (z(\phi\psi)^{s'},w(\psi\phi)^{s'})\\
\iff\; & ((y\psi)(\phi\psi)^{r'},(x\phi)(\psi\phi)^{r'})\sim (z(\phi\psi)^{s'},w(\psi\phi)^{s'})
\end{align*}
which is again an instance of the case for the type VI injective endomorphism $\Psi$ defined above with input $(y\psi,x\phi)$ and $(z,w)$. We can use the same trick to solve the remaining case of deciding if there is a solution $(r,s)\in (2\N+1)\times (2\N+1)$.
We want to  decide whether there are some $r',s'\in \N$ such that $(x,y)\Phi^{2r'+1}\sim (z,w)\Phi^{2s'+1}$.
But 
\begin{align*}
&(x,y)\Phi^{2r'+1}\sim (z,w)\Phi^{2s'+1}\\
\iff\; &  (y(\psi\phi)^{r'}\psi,x(\phi\psi)^{r'}\phi)\sim (w(\psi\phi)^{s'}\psi,z(\phi\psi)^{s'}\phi)\\
\iff\; & ((y\psi)(\phi\psi)^{r'},(x\phi)(\psi\phi)^{r'})\sim ((w\psi)(\phi\psi)^{s'},(z\phi)(\psi\phi)^{s'})
\end{align*}
which is again an instance of the case for the type VI injective endomorphism $\Psi$ defined above with input $(y\psi,x\phi)$ and $(w\psi,z\phi)$. 
 \qed\\

Now we can prove decidability of TCP($F_n\times F_m$) for monomorphisms of $F_n\times F_m$.

\begin{proposition}
There exists an algorithm taking as input two integers $n,m>1$, two elements $(x,y),(z,w)\in F_n\times F_m$ and an injective endomorphism $\Phi\in \Mon(F_n\times F_m)$ that outputs an element $(u,v)\in F_n\times F_m$ such that $(x,y)=(u^{-1},v^{-1})\Phi(z,w)(u,v)$ if such an element exists and outputs \texttt{NO} otherwise.
\end{proposition}
\noindent\textit{Proof.} Suppose that $\Phi$ is a type VI monomorphism. Then, $(x,y)\Phi=(x\phi,y\psi)$ for some $\phi\in \Mon(F_n)$ and $\psi\in \Mon(F_m)$. We have that 
$$(x,y)=(u^{-1},v^{-1})\Phi(z,w)(u,v) \iff \begin{cases}
x=(u^{-1}\phi)zu\\
y=(v^{-1}\psi)wv
\end{cases}.$$
So, there exists an element $(u,v)\in F_n\times F_m$ as desired if and only if TCP($F_n$) answers \texttt{YES} on input $(x,z,\phi)$ and TCP($F_m$) answers \texttt{YES} on input $(y,w,\psi)$. That is decidable by \cite[Lemma 2.2]{[Log22]} or \cite[Theorem 2.4]{[Ven21]}.

Now assume that  $\Phi$ is a monomorphism of type VII. Then $(x,y)\Phi=(y\psi,x\phi)$, for injective homomorphisms $\phi:F_n\to F_m$ and $\psi:F_m\to F_n$. Applying this description of $\Phi$, we obtain that 
\begin{align*}&(x,y)=(u^{-1},v^{-1})\Phi(z,w)(u,v)\iff\begin{cases}
x=(v^{-1}\psi) zu\\
y=(u^{-1}\phi) wv 
\end{cases}\iff
\begin{cases}
u^{-1}=x^{-1}(v^{-1}\psi)z\\
y=(u^{-1}\phi) wv
\end{cases}\\
\iff
&\begin{cases}
u^{-1}=x^{-1}(v^{-1}\psi)z\\
y=(x^{-1}\phi)(v^{-1}\psi\phi)(z\phi) wv
\end{cases}
\iff
\begin{cases}
u=z^{-1}(v\psi)x\\
(x\phi)y=(v^{-1}\psi\phi)(z\phi) wv
\end{cases}.
\end{align*}

So, if $(x,y)$ and $(z,w)$ are $\Phi$-twisted conjugate with conjugator $(u,v)$, then $(x\phi)y$ and $(z\phi) w$ are $\psi\phi$-twisted conjugate with conjugator $v$. Conversely, if $(z\phi) w$ and $(x\phi)y$ are $\psi\phi$-twisted conjugate with some conjugator $v$, then $(x,y)$ and $(z,w)$ are $\Phi$-twisted conjugate with conjugator $(z^{-1}(v\psi)x,v)$. Since it can be decided whether $(z\phi) w$ and $(x\phi)y$ are $\psi\phi$-twisted conjugate, by solving TCP($F_m$), we are done.
\qed\\

By the results in \cite{[Log22]}, the two results above yield the following corollary.
\begin{corollary}\label{conjaschnn}
The conjugacy problem is solvable for ascending HNN-extensions of $F_n\times F_m$, for $n,m>1$.
\end{corollary}
\section{Dynamics at the infinity}

A metric space $(X,d)$ is said to be a \emph{median space} if, for all $x,y,z\in X$, there is some 
unique point $\mu(x,y,z)\in X$, known as the \emph{median} of $x,y,z$, 
such that $d(x,y)=d(x,\mu(x,y,z))+d(\mu(x,y,z),y);$ $d(y,z)=d(y,\mu(x,y,z))+d(\mu(x,y,z),z);$ and $d(z,x)=d(z,\mu(x,y,z))+d(\mu(x,y,z),x)$. We call $\mu:X^3\to X$ the \emph{median operator} of the median space $X$. 

Coarse median spaces were introduced by Bowditch in \cite{[Bow13]}. Following the equivalent definition given in \cite{[NWZ19]}, we say that, given a metric space $X$, a \emph{coarse median on $X$} is a ternary operation $\mu:X^3\to X$ satisfying the following: 

there exists a constant $C\geq 0$ such that, for all $a,b,c,x\in X$, we have that 
\begin{enumerate}
\item $\mu(a,a,b)=a \text{ and } \mu(a,b,c)=\mu(b,c,a)=\mu(b,a,c);$
\item $d(\mu(\mu(a,x,b),x,c),\mu(a,x,\mu(b,x,c)))\leq C;$
\item $d(\mu(a,b,c),\mu(x,b,c))\leq Cd(a,x)+C.$
\end{enumerate} 

Given a group $G$, a \emph{word metric} on $G$ measures the distance of the shortest path in the Cayley graph of $G$ with respect to some set of generators, i.e., for two elements $g,h\in G$, we have that $d(g,h)$ is the length of the shortest word whose letters come from the generating set representing $g^{-1}h$.
Following the definitions in \cite{[Fio21]}, two coarse medians $\mu_1,\mu_2:X^3\to X$ are said to be at \emph{bounded distance} if there exists some constant $C$ such that $d(\mu_1(x,y,z),\mu_2(x,y,z))\leq C$ for all $x,y,z\in X$, and a \emph{coarse median structure} on $X$ is an equivalence class $[\mu]$ of coarse medians pairwise at bounded distance. When $X$ is a metric space and $[\mu]$ is a coarse median structure on $X$, we say that $(X,[\mu])$ is a \emph{coarse median space.} Following Fioravanti's definition in \cite{[Fio21]}, a \emph{coarse median group} is a pair $(G,[\mu])$, where $G$ is a finitely generated group with a word metric $d$ and $[\mu]$ is a $G$-invariant coarse median structure on $G$, meaning that for each $g\in G$, there is a constant $C(g)$ such that $d(g\mu(g_1,g_2,g_3),\mu(gg_1,gg_2,gg_3)\leq C(g)$, for all $g_1,g_2,g_3\in G$. The author in \cite{[Fio21]} also remarks that this definition is stronger than the original definition from \cite{[Bow13]}, that did not require $G$-invariance. Despite being better suited for this work, it is not quasi-isometry-invariant nor commensurability-invariant, unlike Bowditch's version.

Now, a $K$-hyperbolic group is such that there is some $K\geq 0$ for which every geodesic triangle has a $K$-center, i.e., a point that, up to a bounded distance, depends only on the vertices, and is $K$-close to every edge of the triangle. Given three points, the operator that associates the three points to the $K$-center of a geodesic triangle they define is coarse median. In fact, by \cite[Theorem 4.2]{[NWZ19]} it is the only coarse-median structure that we can endow $X$ with.
Given two groups $G_1$ and $G_2$ endowed with coarse median operators $\mu_1$ and $\mu_2$, then it is easy to check that the operator $\mu:(G_1\times G_2)^3\to G_1\times G_2$ defined by $\mu((x_1,x_2),(y_1,y_2),(z_1,z_2))=(\mu_1(x_1,y_1,z_1),\mu_2(x_2,y_2,z_2))$ is coarse median. We refer to $\mu$ as the \emph{product coarse median operator} and this is the coarse median we will endow $F_n\times F_m$ with.

Given a coarse median group $(G,[\mu])$, an automorphism $\varphi\in \Aut(G)$ is said to be \emph{coarse-median preserving} if it fixes $[\mu]$, i.e., when there is some constant $C\geq 0$ such that for all $g_i'$s, we have that $$d(\mu(g_1,g_2,g_3)\varphi,\mu(g_1\varphi,g_2\varphi,g_3\varphi))\leq C,$$ with respect to some word metric $d$. This can naturally be defined for general endomorphisms, not necessarily bijective, as done in \cite{[Car22a],[Car22c]}, and even to homomorphisms between different groups: we say that a homomorphism $\varphi:(G_1,\mu_1)\to (G_2,\mu_2)$ is coarse-median preserving if there is some $C\geq 0$ such that for all $x,y,z\in G_1$,
$$d(\mu_2(x\varphi,y\varphi,z\varphi),\mu_1(x,y,z)\varphi)\leq C$$
for some word metric $d$ of $G_2$.

Let $G=\langle A\rangle$ be a hyperbolic group. Given $g,h,p\in G$, we define the \emph{Gromov product of $g$ and $h$} taking $p$ as basepoint by
$$(g|h)_p^A=\frac 12 (d_A(p,g)+d_A(p,h)-d_A(g,h)).$$ We will  write $(g|h)$ to denote $(g|h)_1^A$, when the generating set is clear. Notice that, in the free group case, we have that $(g|h)=|g\wedge h|,$ where $u\wedge v$ denotes the longest common prefix between $u$ and $v$.

 Recall the notation from Section 2, namely the definition of sets $X, \, Y,\,Z$ and $W$.

It is shown in \cite[Theorem 4.11]{[Car22c]} that an endomorphism $\phi$ of a hyperbolic group $G$ is coarse-median preserving if and only if the bounded reduction property holds for $\phi$. In particular, if $G$ is free, $\phi$ is coarse-median preserving if and only if it is either trivial or injective. Following the same strategy as in \cite[Theorem 4.11]{[Car22c]}, we can see that the result holds for homomorphisms between free groups of different ranks.

\begin{proposition}\label{cmphomo}
Let $\varphi:F_n\to F_m$ be a homomorphism. Then $\varphi$ is coarse-median preserving if and only if it is either trivial or injective.
\end{proposition}
\noindent\textit{Proof.} Suppose that $\varphi$ is not trivial nor injective. Then there are some nontrivial $x,y\in F_n$ such that $x\varphi=1$ and $y\varphi\neq 1$. Eventually conjugating $x$ by a suitable element, we can assume that $y^Cxy^{-C}$ is reduced for all $C\in \N$. Indeed, conjugating a kernel element by any element yields a kernel element, and if $c$ is the last letter of $y$, then $c^{|x|}xc^{-|x|}$ is an example of an element with the desired property (notice that we may assume that $x$ is not itself a power of $c$, as that would imply that $c\varphi$ is $1$, in which case we could simply remove all occurrences of $c$ and $c^{-1}$ from $y$).
 So, 
for all $C\in \N$, we have that 
\begin{align*}
&d(\mu_1(x,y^Cxy^{-C},y^C)\varphi,\mu_2(x\varphi,(y^Cxy^{-C})\varphi,y^C\varphi))\\
=\; &d(y^C\varphi,\mu_2(1,1,y^C\varphi))\\
=\; & d(1,y^C\varphi)\\
\geq\; &C.
\end{align*}
Since $C$ can be arbitrarily large, $\varphi$ is not coarse-median preserving.

Clearly, if $\varphi$ is trivial, it is coarse-median preserving, so the only thing left to prove is that injectivity implies coarse-median preservation. Suppose that $\varphi$ is injective. By \cite[Proposition 5.3]{[Car22c]}, for all $M\geq 0$, there is some $N\geq 0$ such that $$(u|v)\leq M \implies (u\varphi |v\varphi)\leq N,$$ for all $u,v\in F_n$. Following the proofs of \cite[Proposition 4.2 and Theorem 4.6]{[Car22c]} step by step, we get that there is some $N\in \N$ such that, for all $x,y\in F_n$, the image of the geodesic $[x,y]$ is at a Hausdorff distance at most $N$ from the geodesic $[x\varphi,y\varphi]$. Let $x,y,z\in F_n$. Then $\mu_1(x,y,z)$ belongs to $[x,y]$, $[y,z]$ and $[x,z]$, and so $\mu_1(x,y,z)\varphi$ is at a distance at most $N$ from $[x\varphi,y\varphi]$, $[y\varphi,z\varphi]$ and $[x\varphi,z\varphi]$, thus $\mu_1(x,y,z)\varphi$ is an $N$-center of the triangle $[[x\varphi,y\varphi,z\varphi]]$. Since $\mu_2(x\varphi,y\varphi,z\varphi)$ is also an $N$-center of $[[x\varphi,y\varphi,z\varphi]]$, by  \cite[Lemma 3.1.5]{[Bow91]}, it is at a bounded distance from $\mu_1(x,y,z)\varphi$ and this bound depends only on $N$.

Therefore, $\varphi$ is coarse-median preserving.
\qed\\

Now we establish an equivalence between preservation of a coarse median and uniform continuity in $F_n\times F_m$, similar to what happens in free and free-abelian times free groups.
\begin{theorem}\label{uccmp}
Coarse-median preserving endomorphisms of $F_n\times F_m$ are precisely the uniformly continuous ones.
\end{theorem}
\noindent\textit{Proof.} Let $\Phi\in \End(F_n\times F_m)$ be a coarse-median preserving endomorphism  with constant $C$ and  recall the notation from Section 2, namely the definition of sets $X, \, Y,\,Z$ and $W$.

We start by proving that $1\in X\implies X=\{1\}$. Suppose there are some $i,j\in [n]$ such that $(a_i,1)\Phi=(1,y_i)$ and $(a_j,1)\Phi= (x_j,y_j)$ with $x_j\neq 1$. Then, letting $D$ be a number greater than $C$,
\begin{align*}
&d\left(\mu\left((a_j^{-D}a_i,1)\Phi,(a_j^{-D}a_ia_j^D,1)\Phi,(a_i,1)\Phi\right),\mu\left((a_j^{-D}a_i,1),(a_j^{-D}a_ia_j^D,1),(a_i,1)\right)\Phi\right)\\
=\;&d\left(\mu\left((x_j^{-D},y_j^{-D}y_i),(1,y_j^{-D}y_iy_j^D),(1,y_i)\right),\left(\mu_1(a_j^{-D}a_i,a_j^{-D}a_ia_j^D,a_i),\mu_2(1,1,1)\right)\Phi\right)\\
=\;&d\left(\left(\mu_1(x_j^{-D},1,1),\mu_2(y_j^{-D}y_i,y_j^{-D}y_iy_j^D,y_i)\right),(a_j^{-D}a_i,1)\Phi\right)\\
\geq \;& d(\mu_1(x_j^{-D},1,1),x_j^{-D})\\
=\;&d(1,x_j^{-D})\\
>\;&C
\end{align*}
which contradicts the fact that $\Phi$ is coarse-median preserving with constant $C$. So, $1\in X\implies X=\{1\}$. 
The same holds for $Y, \,Z$ and $W$, and that can be seen analogously. Notice that in the second equality above, we have that 
$\mu_1(a_j^{-D}a_i,a_j^{-D}a_ia_j^D,a_i)=a_j^{-D}a_i$ since the elements $a_i$ are distinct basis elements, hence all words appearing in the expression are freely reduced. 

We now prove that we cannot have $Y$ and $W$ to be nontrivial simultaneously, which shows that $\Phi$ cannot be of type I, II or V. Suppose that $Y\neq\{1\}$ and $W\neq \{1\}$. In view of the above, $1\not\in Y\cup W$. We know that there is some word 
$v\in F_m$ and nonzero exponents $q_i$ and $s_j$ such that $y_i=v^{q_i}$ and $w_j=v^{s_j}$.
For $D\in \N$, we have that 
\begin{align*}
&\mu\left((a_i^{Ds_j},b_j^{-Dq_i})\Phi,(a_i^{\sgn(s_j)+2Ds_j},b_j^{-Dq_i})\Phi,(a_i^{2Ds_j},b_j^{-2Dq_i})\Phi\right)\\
=\;&\mu\left((x_i^{Ds_j}z_j^{-Dq_i},1),(x_i^{2Ds_j+\sgn(s_j)}z_j^{-Dq_i},v^{q_i(\sgn(s_j)+Ds_j)}),(x_i^{2Ds_j}z_j^{-2Dq_i},1)\right)\\
=\; & \left(\mu_1(x_i^{Ds_j}z_j^{-Dq_i},x_i^{2Ds_j\sgn(s_j)}z_j^{-Dq_i},x_i^{2Ds_j}z_j^{-2Dq_i}),\mu_2(1,v^{q_i(\sgn(s_j)+Ds_j)},1)\right)
\end{align*}
and 
\begin{align*}
&\mu\left((a_i^{Ds_j},b_j^{-Dq_i}),(a_i^{\sgn(s_j)+2Ds_j},b_j^{-Dq_i}),(a_i^{2Ds_j},b_j^{-2Dq_i})\right)\Phi\\
=\;&\left(\mu_1(a_i^{Ds_j},a_i^{\sgn(s_j)+2Ds_j},a_i^{2Ds_j}),\mu_2(b_j^{-Dq_i},b_j^{-Dq_i},b_j^{-2Dq_i})\right)\Phi\\
=\; &(a_i^{2Ds_j},b_j^{-Dq_i})\Phi.
\end{align*}

Then, for $D>C$, the distance between $$\mu\left((a_i^{Ds_j},b_j^{-Dq_i})\Phi,(a_i^{\sgn(s_j)+2Ds_j},b_j^{-Dq_i})\Phi,(a_i^{2Ds_j},b_j^{-2Dq_i})\Phi\right)$$
and
$$\left(\mu((a_i^{Ds_j},b_j^{-Dq_i}),(a_i^{\sgn(s_j)+2Ds_j},b_j^{-Dq_i}),(a_i^{2Ds_j},b_j^{-2Dq_i})\right)\Phi$$
is greater than or  equal to $$d(\mu_2(1,v^{q_i(\sgn(s_j)+Ds_j)},1),v^{Ds_jq_i})=d(1,v^{Ds_jq_i})\geq C,$$
which contradicts the fact that $\Phi$ is coarse-median preserving with constant $C$. The same argument yields that we cannot have both $X$ and $Z$ nontrivial, so type III is also dealt with.

Since coarse-median and uniform continuity coincide for free groups, it follows from \cite[Lemma 3.2]{[Car22a]} that an endomorphism of type VI is coarse-median preserving if and only if it is uniformly continuous.

If $\Phi$ is a type IV endomorphism defined as $(x,y)\Phi=(y\phi,y\psi)$, then, for all $(x_1,x_2),(y_1,y_2),(z_1,z_2)\in F_n\times F_m$ and $C>0$, we have that 
\begin{align*}
&d(\mu((x_1,x_2),(y_1,y_2),(z_1,z_2))\Phi,\mu((x_1,x_2)\Phi, (y_1,y_2)\Phi,(z_1,z_2)\Phi))\leq C\\
\iff\; &d((\mu_1(x_1,y_1,z_1),\mu_2(x_2,y_2,z_2))\Phi,\mu((x_2\phi,x_2\psi), (y_2\phi,y_2\psi),(z_2\phi,z_2\psi)))\leq C\\
\iff\; & d((\mu_2(x_2,y_2,z_2)\phi,\mu_2(x_2,y_2,z_2)\psi),(\mu_1(x_2\phi,y_2\phi,z_2\phi),\mu_2(x_2\psi,y_2\psi,z_2\psi)))\leq C\\
\iff\; & \begin{cases}
d((\mu_2(x_2,y_2,z_2)\phi,\mu_1(x_2\phi,y_2\phi,z_2\phi))\leq C\\
d(\mu_2(x_2,y_2,z_2)\psi,\mu_2(x_2\psi,y_2\psi,z_2\psi))\leq C
\end{cases}.
\end{align*}
Hence $\Phi$ is coarse-median preserving if and only if both $\phi$ and $\psi$ are. Since, in view of Proposition \ref{cmphomo} this is equivalent to having  injective or trivial $\phi$ and $\psi$. These are precisely the uniformly continuous type IV endomorphisms of $F_n\times F_m$ (see \cite[Proposition 6.2]{[Car23b]}).

Finally, if $\Phi$ is a type VII endomorphism defined as $(x,y)\Phi=(y\psi,x\phi)$, then, 
proceeding as in the type IV case, we obtain that $\Phi$ is coarse-median preserving if and only if both $\phi$ and $\psi$ are. Since, in view of Proposition \ref{cmphomo} this is equivalent to having  injective or trivial $\phi$ and $\psi$. These are precisely the uniformly continuous type VII endomorphisms of $F_n\times F_m$ (see \cite[Proposition 6.2]{[Car23b]}).
\qed\\

We will now denote by $d_1$ the prefix metric on $F_n$, by $d_2$ the prefix metric on $F_m$ and by $d$ the product metric obtained by taking $d_1$ and $d_2$ in the factors. We can now prove the main result of the section.

\begin{theorem}\label{recper}
The continuous extension of every automorphism of $F_n\times F_m$ to the completion obtained by taking the prefix metric in each factor has asymptotically periodic dynamics.
\end{theorem}
\noindent\textit{Proof.} We will prove that  a point $(z,w)\in \widehat{F_n\times F_m}$ that belongs to the $\omega$-limit of some point $(x,y)$ must be periodic, which implies that the $\omega$-limit is itself the (periodic) orbit of $(z,w)$, and so that the dynamics is asymptotically periodic.

Let $\Phi\in \Aut(F_n\times F_m)$ be a type VI automorphism of $F_n\times F_m$. Then $\Phi$ is of the form $(x,y)\mapsto (x\phi,y\psi)$, for some $\phi\in \Aut(F_n)$ and $\psi\in \Aut(F_m)$. By \cite[Proposition 6.2]{[Car23b]}, $\Phi$ is uniformly continuous with respect to the metric $d$ and, by uniqueness of the extension $\widehat\Phi$ is given by $(x,y)\mapsto (x\widehat\phi, y\widehat\psi)$. Let $(z,w)\in \omega((x,y),\widehat\Phi)$. There is a sequence  $(n_i)_{i\in \N}$ such that $(x,y)\widehat\Phi^{n_i}\to (z,w)$. This means that 
  $$\left( x\widehat\phi^{n_i},y\widehat\psi^{n_i}\right)\to (z,w),$$ 
  and so 
 $ x\widehat\phi^{n_i}\to z$ and $y\widehat\psi^{n_i}\to w$, which implies that $z\in \omega(x,\widehat\phi)$ and $w\in \omega(y,\widehat\psi)$. Since the dynamics of free group automorphisms is asymptotically periodic by Theorem \ref{levittlustig}, then $x$ must be $\widehat\phi$-periodic and $y$ must be $\widehat\psi$-periodic, and so $(x,y)$ is $\widehat\Phi$-periodic.

Now, suppose that $n=m$ and let $\Phi\in \Aut(F_n\times F_m)$ be a type VII automorphism of $F_n\times F_n$.  Then $\Phi$ is of the form $(x,y)\mapsto (y\psi,x\phi)$, for some $\phi,\psi\in \Aut(F_n)$. Again, by \cite[Proposition 6.2]{[Car23b]}, $\Phi$ is uniformly continuous with respect to the metric $d$ and, by uniqueness of the extension $\widehat\Phi$ is given by $(x,y)\mapsto (y\widehat\psi, x\widehat\phi)$.
Moreover, $\widehat \phi$ and $\widehat\psi$ are continuous maps between compact spaces, hence uniformly continuous and since both $\phi$ and $\psi$ are automorphisms, then $\phi^{-1}$ and $\psi^{-1}$ are also uniformly continuous and extend to continuous mappings $\widehat{\phi^{-1}}$ and $\widehat{\psi^{-1}}$. By uniqueness of extensions, we have that $\widehat\phi$ and $\widehat{\phi^{-1}}$ (resp. $\widehat\psi$ and $\widehat{\psi^{-1}}$) are mutually inverse uniformly continuous mappings.

Let $(z,w)\in \omega((x,y),\widehat\Phi)$. Then, there is a sequence $(n_i)_{i\in \N}$ such that $(x,y)\widehat\Phi^{n_i}\to (z,w)$. If $(n_i)$ has infinitely many even numbers, then consider the subsequence $(n_{k_i})_{i\in \N}$ of even numbers such that  $(x,y)\widehat\Phi^{n_{k_i}}\to (z,w)$. This means that 
  $$\left( x(\widehat\phi\widehat\psi)^{\frac{n_{k_i}}{2}},y(\widehat\psi\widehat\phi)^{\frac{n_{k_i}}{2}}\right)\to (z,w),$$ and so 
 $x(\widehat\phi\widehat\psi)^{\frac{n_{k_i}}{2}}\to z$ and $y(\widehat\psi\widehat\phi)^{\frac{n_{k_i}}{2}}\to w$, which implies that $z\in \omega(x,\widehat\phi\widehat\psi)$ and that $w\in \omega(y,\widehat\psi\widehat\phi)$. By Theorem \ref{levittlustig}, we deduce that 
$z$ is $(\widehat\phi\widehat\psi)$-periodic and   $w$ is $(\widehat\psi\widehat\phi)$-periodic, that is, that 
there are $p_1,p_2\in \N$ such that $z(\widehat\phi\widehat\psi)^{p_1}=z$ and $w(\widehat\psi\widehat\phi)^{p_2}=w$. Hence, $$(z,w)\widehat\Phi^{2p_1p_2}=(z(\widehat\phi\widehat\psi)^{p_1p_2},w(\widehat\psi\widehat\phi)^{p_1p_2})=(z,w),$$
that is $(z,w)$ is $\widehat\Phi$-periodic.

  If $(n_i)$ has only finitely many even numbers, then it must have infinitely many odd numbers and we consider the subsequence $(n_{k_i})_{i\in \N}$ of odd numbers such that  $(x,y)\widehat\Phi^{n_{k_i}}\to (z,w)$.
  This means that 
  $$\left( y(\widehat\psi\widehat\phi)^{\frac{n_{k_i}-1}{2}}\widehat\psi,x(\widehat\phi\widehat\psi)^{\frac{n_{k_i}-1}{2}}\widehat\phi\right)\to (z,w),$$ 
  and so $y(\widehat\psi\widehat\phi)^{\frac{n_{k_i}-1}{2}}\widehat\psi\to z$ and $x(\widehat\phi\widehat\psi)^{\frac{n_{k_i}-1}{2}}\widehat\phi\to w$. By uniform continuity of $\widehat\phi$ and $\widehat\psi$, it follows that 
   $y(\widehat\psi\widehat\phi)^{\frac{n_{k_i}-1}{2}}\widehat\psi\widehat\phi\to z\widehat\phi$ and $x(\widehat\phi\widehat\psi)^{\frac{n_{k_i}-1}{2}}\widehat\phi\widehat\psi\to w\widehat\psi$. Hence, 
   $z\widehat\phi\in \omega(y,\widehat\psi\widehat\phi)$ and $w\widehat\psi\in \omega(x,\widehat\phi\widehat\psi)$, which, by Theorem \ref{levittlustig}, implies that 
$z\widehat\phi$ is $(\widehat\psi\widehat\phi)$-periodic and   $w\widehat\psi$ is $(\widehat\phi\widehat\psi)$-periodic, that is, that there are $p_1,p_2\in \N$ such that $z\widehat\phi(\widehat\psi\widehat\phi)^{p_1}=z\widehat\phi$ and $w\widehat\psi(\widehat\phi\widehat\psi)^{p_2}=w\widehat\psi$. 
Therefore,  $$z(\widehat\phi\widehat\psi)^{p_1}=z\widehat\phi(\widehat\psi\widehat\phi)^{p_1}\widehat{\phi^{-1}}=z\widehat\phi\widehat{\phi^{-1}}=z$$ and $$w(\widehat\psi\widehat\phi)^{p_2}=w\widehat\psi(\widehat\phi\widehat\psi)^{p_2}\widehat{\psi^{-1}}=w\widehat\psi\widehat{\psi^{-1}}=w.$$
Hence, $$(z,w)\widehat{\Phi}^{2p_1p_2}=(z(\widehat\phi\widehat\psi)^{p_1p_2},w(\widehat\psi\widehat\phi)^{p_1p_2})=(z,w),$$
and $(z,w)$ is $\widehat\Phi$-periodic. 
 \qed

\section*{Acknowledgements}
The author thanks Rémi Coulon, Damien Gaboriau, and Yassine Guerch for pointing out a flaw in a previous version of the paper. The author is also indebted to the anonymous referee for their report, which confirmed this issue and offered further comments and suggestions that helped improve the overall quality of the paper.
This work is funded by national funds through the FCT - Fundação para a Ciência e a Tecnologia, I.P., under the scope of the projects UIDB/00297/2020 and UIDP/00297/2020 (Center for Mathematics and Applications).

\bibliographystyle{plain}
\bibliography{Bibliografia}

 \end{document}